\documentclass{amsart}
\input{amssymb.sty}

\newcommand{\PSbox}[3]{\mbox{\rule{0in}{#3}\hspace{#2}\includegraphics{#1}}}

\numberwithin{equation}{section}
\newtheorem{thm}{Theorem}
\newtheorem{lemma}[thm]{Lemma}

\newtheorem{prop}[thm]{Proposition}

\newcommand{\R}{{\mathbb R}}
\newcommand{\C}{{\mathbb C}}
\newcommand{\Z}{{\mathbb Z}}
\newcommand{\Q}{{\mathbb Q}}
\newcommand{\E}{{\bf E}}
\renewcommand{\P}{{\bf P}}

\title{Local statistics of lattice dimers}
\author{Richard Kenyon}
\address{CNRS UMR 8628, Laboratoire de Math\'ematiques,
Universit\'e Paris-Sud, 91405 Orsay, France.}
\begin{document}

\begin{abstract}
We show how to compute the probability of any given local configuration
in a random tiling of the plane with dominos.
That is, we explicitly compute
the measures of cylinder sets for the measure of maximal entropy $\mu$ 
on the space of tilings of the plane with dominos.

We construct a measure $\nu$ on the set of lozenge tilings of the plane,
show that its entropy is the topological entropy,
and compute explicitly the $\nu$-measures of cylinder sets.

As applications of these results,
we prove that the translation action is strongly mixing for $\mu$ and $\nu$,
and compute the rate of convergence to mixing 
(the correlation between distant events).
For the measure $\nu$ we compute the variance of the height function.
\medskip

\noindent{\bf Resum\'e.} Soit $\mu$ la mesure d'entropie maximale
sur l'espace $X$ des pavages du plan par des dominos. On calcule
explicitement la mesure des sous-ensembles cylindriques de $X$.
De m\^eme,
on construit une mesure $\nu$ d'entropie maximale 
sur l'espace $X'$ des pavages du plan
par losanges, et on calcule explicitement 
la mesure des sous-ensembles cylindriques.

Comme application on calcule, pour $\mu$ et $\nu$, les correlations
d'\'evenements distants, ainsi que la $\nu$-variance
de la fonction ``hauteur" sur $X'$.
\end{abstract}

\maketitle
\section{Introduction}
A {\bf domino} is
a $2$ by $1$ or a $1$ by $2$ 
rectangle, whose vertices have integer coordinates in the plane.
Let $X$ be the space of all tilings of the plane with dominos.
Then $X$ has a natural topology, where 
two tilings are close if they agree
on a large neighborhood of the origin. In this topology, $X$ is compact,
and $\Z^2$ acts 
continuously on $X$ by translations.

Burton and Pemantle \cite{BP} proved that there is a unique invariant 
measure $\mu$ of maximal entropy for this action.
This measure has the following property:
let $T$ be any tiling of a finite region $R$, and $U_T$ the set of 
tilings in $X$ extending $T$ (we call $U_T$ 
{\bf the cylinder set generated by $T$}).
Then $\mu(U_T)$ only depends on $R$ and not on the choice $T$
of a tiling of $R$.

In this paper we compute $\mu(U_T)$ explicitly, for any finite
tiling $T$.  We will prove in section \ref{dominopf} that:

\begin{thm}\label{domc}
There is a function $P\colon \Z^2\to\C$ with the following property.
Let $E$ be any finite set of disjoint dominos,
covering ``black" squares $b_1,\ldots,b_k$ and ``white" squares
$w_1,\ldots,w_k$. Then $\mu(U_E)$ is the absolute value of the determinant
of the $k\times k$ matrix $M=(m_{ij})$, where $m_{ij}=P(b_i-w_j)$.
(Here $b_i-w_j$ is the translation vector from square $w_j$ to square $b_i$.)
\end{thm}

We call $P$ the {\bf coupling function} for $\mu$.
It can be computed explicitly (see Theorem \ref{Pdoms} below),
and the values of $P$ are in $\Q\oplus\frac1\pi\Q$.

A similar result holds for tilings of the plane with lozenges.
A {\bf lozenge} is a rhombus with side $1$, smaller angle $\pi/3$,
and vertices in the lattice $\Z[e^{2\pi i/3}]$.
We define a measure $\nu$ on the space $X'$ of lozenge tilings,
which is the limit of uniform measures on periodic tilings,
and show that its entropy equals the topological entropy of the $\Z^2$-action.
(For a background on measure-theoretic and topological entropy see
\cite{Petersen}.)
We conjecture that $\nu$ is the {\it unique} measure of maximal entropy.
We compute the coupling function for $\nu$ in Theorem \ref{lozcouple}:
it takes values in $\Q\oplus\frac{\sqrt{3}}{\pi}\Q$.
\medskip

We give several applications of these two results.
First, it is known that the translation-action of $\Z^2$ on dominos or
lozenges is topologically mixing (if two events are far enough apart, one can
find a single tiling containing both of them). 
We show here that it is {\bf strongly mixing} for the measures $\mu$ and $\nu$,
that is, events that are far apart are almost independent.
In particular we show that the convergence rate to independence 
is at least quadratic in the 
distance, that is,  for any two finite tilings $T_1$
and $T_2$, and $v\in\Z^2$,
$$\mu(U_{T_1}\cap U_{(v+T_2)})= \mu(U_{T_1})\mu(U_{T_2})+O(\frac1{|v|^2}),$$ 
and similarly for $\nu$;
see Theorem \ref{mixing}. There is a corresponding lower bound
when $T_1,T_2$ are single tiles.
Correlations for single tiles 
have been previously computed by Fisher and Stephenson \cite{FS} and
Stephenson \cite{St}.
 
Another application 
is to the computation of the variance of the height function for $\nu$.
The {\bf height function} of a lozenge tiling was introduced by
Bl\"ote and Hilhorst \cite{BH}.
We compute for $\nu$ the variance of the height function in the $y$-direction 
(Theorem \ref{var}), that is, 
the variance in height between two vertices separated by distance
$n$ in the $y$-direction. 
This variance is related to the expected number of ``contour lines" separating
the two points (Proposition \ref{contours}).
The variance in other directions can also be computed with the same methods,
although we could not obtain a closed form.
On the other hand, domino tilings also have a height function, but we were not
able to compute its variance using the methods here.

A third application which we will not discuss here
is to the computation of entropy $\mbox{ent}(s,t)$ of tilings
whose height function has a fixed average slope, that is,
the set of tilings of the plane 
whose height function $ht(\cdot,\cdot)$ satisfies
$ht(x,y)=ht(0,0)+sx+ty+o(|x|,|y|)$. This computation is used in \cite{CKP}
to give a method of counting the approximate number of tilings of any region,
by maximizing an entropy functional over the region.

The results of this paper are very general.
Using a result of Kasteleyn \cite{Kast2}, we can
compute cylinder set measures (for certain
measures of maximal entropy) for the space of perfect matchings
on any periodic planar graph. Here ``periodic" means
the graph is finite modulo an action of $\Z^2$ by translations.
In particular Hurst and Green \cite{HG}
showed that the {\bf Ising model} on a planar graph
with nearest-neighbor interactions
can be represented as a matching problem on another planar graph
(actually they showed this in an equivalent form).
Our method
therefore allows one to compute the local probabilities for the Ising model.
Indeed, methods similar to those in this paper
were used in Montroll, Potts and Ward \cite{MPW} to compute
long-range correlations for the Ising model.
\medskip

The coupling function which we define is 
analogous to the Green's function for a resistor network
(see \cite{BP}). Indeed, subdeterminants of the Green's function matrix
give probabilities of finding subsets of edges in a random spanning tree.
There is a well-known connection between domino tilings and spanning trees:
Burton and Pemantle \cite{BP} have found,
via enumeration of spanning trees, the measures of certain
cylinder sets for dominos. Our method, which is different,
gives the measure of all cylinder
sets and is adaptable to planar graphs where the spanning tree
connection is not known to hold.

The structure of this paper is as follows.
Section \ref{bkgd} discusses background: we show how to compute the number
of tilings as a determinant, and we define $\mu,\nu$ and the height
function. In section \ref{cyls} we show how to compute the cylinder set measures
for finite planar graphs. In section \ref{lozs} we compute explicitly
the coupling function for $\nu$ and give the applications.
We chose in this section 
to concentrate on the case of lozenges rather than dominos since lozenges
seem to be less prominent in the mathematics literature. 
Section \ref{dominopf} restates the results of section \ref{lozs} for dominos.

\section{Background}
\label{bkgd}
\subsection{Matchings on finite planar graphs}
Let $G$ be a graph.
Two edges of $G$ are {\bf disjoint} if they have no vertices in common.
A {\bf perfect matching} of $G$ is a set of disjoint edges
which contains all the vertices. If $G$ is finite, connected and bipartite,
we say $G$ is {\bf balanced} if there are the same number of vertices in
each half of the bipartition.

Let $\Z^2$ denote the graph whose vertices are $\Z^2$ and whose edges
join every pair of vertices of distance $1$. 
It is easy to see that domino tilings of the plane correspond bijectively
to perfect matchings of $\Z^2$.
The graph $\Z^2$ is bipartite: color the vertex $(x,y)$ black if $x+y\equiv
0\bmod 2$, otherwise color it white. Note that
then a domino ``covers" exactly one vertex
of each color.

\begin{figure}[htbp]
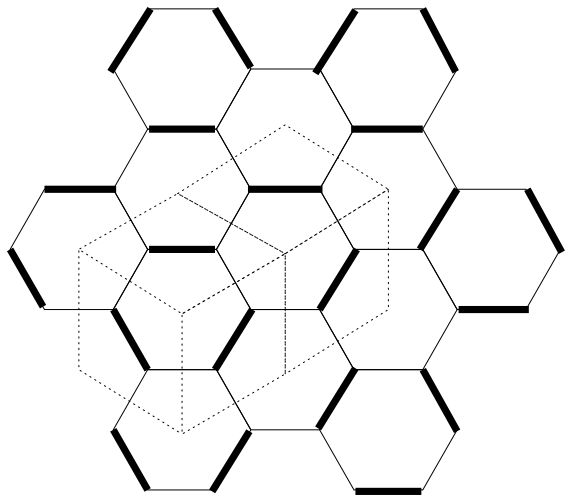

\PSbox{hexgraph.eps}{1in}{3in}
\caption{\label{hexgraph}A matching on the graph $H$ and some of the
corresponding lozenges (dotted lines).}
\end{figure}

Similarly, lozenge tilings correspond bijectively to perfect matchings of a 
particular bipartite planar graph $H$, the ``honeycomb'' grid. 
The graph $H$ is the $1$-skeleton of the
periodic tiling of the plane with regular hexagons (see Figure \ref{hexgraph}).
Rather than use lozenges as defined in the introduction, we adopt the
following coordinates.
Let $\hat x=\frac32-\frac{i\sqrt{3}}2$ and 
$\hat y=\frac32+\frac{i\sqrt{3}}2$.
The set of vertices of $H$ is $V=V_0\cup V_1$, where 
$V_0=\{a\hat{x}+b\hat{y}~\mid~a,b\in\Z\}$  are the ``black" vertices,
and 
$V_1=1+V_0$ are the ``white" vertices.
There is an edge between every pair of vertices at unit distance.
We will denote the vertex $a\hat{x}+b\hat{y}$ by the triple $(a,b,0)$,
and a vertex $a\hat{x}+b\hat{y}+1$ by the triple $(a,b,1)$. 
Hopefully denoting vertices with ordered triples will not confuse the reader.
Note also that a lozenge has side $\sqrt{3}$ in these coordinates.

Kasteleyn \cite{Kast1} and Temperley and Fisher \cite{TF} independently
computed the number of perfect matchings of an $m\times n$ grid 
(in other words, the number of domino tilings of an $m\times n$ rectangle).
Kasteleyn in \cite{Kast1} also computed the number of matchings on a 
toroidal graph in the case of dominos. Later, 
in \cite{Kast2}, he
showed how to compute the number of perfect matchings of general
planar graphs (thereby including the case of lozenges).
His method is to count tilings using Pfaffians. 
We outline here an adaptation of his proof for the case of a 
simply-connected subgraph of $H$.

We say that a subgraph of $H$ is {\bf simply connected} if it
is the $1$-skeleton of a simply connected union of 
``basic hexagons'' of $H$. The following is a special case of the results
of \cite{Kast2}:
\begin{thm}
Let $H'$ be a finite simply connected subgraph of $H$.
The number of perfect matchings  of $H'$
equals $\sqrt{|\det(A)|}$, where $A$ is the adjacency matrix of $H'$.
\end{thm}

Here is a sketch of the proof:
since the graph $H'=(V',E')$ is bipartite, we have $V'=V'_0\cup V'_1$
(where $V_0'\subset V_0, V_1'\subset V_1$)
and $E'\subset V'_0\times V'_1$.
Order the vertices $V'$ so that those in $V'_0$ comes first.
We then can write
\begin{equation}\label{B}
A=\left(\begin{array}{cc}0&B\\B^t&0\end{array}\right),
\end{equation}
where $B:\R^{|V'_0|}\to\R^{|V'_1|}$.
We can assume that $V'_0$ and $V'_1$
have the same cardinality $n$ (that is, $H'$ is balanced),
for otherwise there are no perfect matchings (and $\det(A)$ is clearly $0$).
Then the determinant of the square matrix $B$ is
\begin{equation}\label{detB}\det B=\sum_{\sigma}
\mbox{sgn}(\sigma)a_{1\sigma(1)}\cdots a_{n\sigma(n)}
\end{equation}
where $A=(a_{ij})$ and the sum is over bijections
$\sigma\colon V'_0\to V'_1$.
Each {\em nonzero} term in the sum corresponds to a matching.
One need only check (and this is the most interesting part,
although we won't do it here) that each term has the same sign.
Thus $|\det(B)|=\sqrt{|\det(A)|}={\rm~\#~of~matchings}$.

This same method also works for dominos,
except that one must modify the adjacency matrix
to make the signs in (\ref{detB}) work out right. One way to modify it
is to put weight $i=\sqrt{-1}$ on vertical edges \cite{Wu}.

It is slightly more complicated to count matchings on a 
graph embedded on a torus.  Kasteleyn showed that, for a graph on a torus, 
one can count matchings using a sum of $4$ Pfaffians (see Lemma \ref{BBBB} below).

\subsection{The measure of maximal entropy}
An important issue not dealt with by Kasteleyn is the appropriate sense
in which random tilings of large finite regions in the plane resemble
random tilings of the whole plane. 

Indeed, as the work of \cite{EKLP} shows, the local configuration
of the boundary 
of a region has a drastic effect on the number of perfect matchings,
and hence one must take care when computing the 
measure $\mu$ as a limit of measures on tilings of finite regions.

A {\bf tiling with free boundary conditions}, or simply {\bf free tiling},
of a region $R$ is a tiling of $R$, where the tiles are allowed
to protrude beyond the boundary of $R$.
The {\bf topological entropy} of domino tilings is by definition
$$H_{top}=\lim_{k\to\infty}\frac{1}{k^2}\log|S_k|,$$
where $S_k$ is the set of free tilings of a $k\times k$ square.
Here the square regions  $S_k$
can be replaced by any sequence of ``sufficiently nice" 
regions, for example, convex regions containing squares of size
tending to infinity.  In this case
one divides by the area of the $k$-th region, not by $k^2$.
The same definition also applies to lozenge tilings (replace the $k\times k$
square with for example the set of vertices in a $k\times k$ square
centered at the origin).

The {\bf entropy} $H(\mu')$ of a translation-invariant 
measure $\mu'$ on domino tilings of the plane
can be defined by
$$H(\mu')= \lim_{k\to\infty}\frac{1}{k^2}\sum_{T\in S_k} 
-\mu'(U_T)\log\mu'(U_T).$$
We have $H(\mu')\leq H_{top}$ for any translation-invariant
measure $\mu'$ on $X$, since, as the reader may readily show,
for a finite set of real numbers $\{a_i\}$ satisfying $a_i>0$ and $\sum a_i=1$, 
the quantity $-\sum a_i\log a_i$ is maximized when the
$a_i$ are equal.

For any measure $\mu'$ defined on the finite set $X(S)$ 
of tilings of a finite region $S$, the entropy (per site) of $\mu'$ is
$$H(\mu')=-\frac1{|S|}\sum_{s\in X(S)} \mu'(s)\log\mu'(s),$$
where $\mu'(s)$ is the probability of the tiling $s$.

\subsubsection{The domino measure}
Burton and Pemantle showed:
\begin{thm}[\cite{BP}]\label{unique}
There is a unique measure of maximal entropy $\mu$ for domino
tilings of the plane. 
The entropy of $\mu$ equals the topological entropy $H_{top}(X)$.
\end{thm}
Let $\Z^2_{m,n}$ be the quotient of $\Z^2$ by the lattice
of translations generated by $(2m,0)$ and $(0,2n)$. Then
$\Z^2_{m,n}$ is a graph (on a torus) with $4mn$ vertices.
Let $\mu_{m,n}$ denote the uniform measure on perfect matchings of
$\Z^2_{m,n}$.

Let $\mu^*$ be a weak limit of $\mu_{m,n}$ as $m,n\to\infty$.
We claim that $\mu^*=\mu$.
To see this, Kasteleyn \cite{Kast1} showed that the entropy of $\mu_{m,n}$
converges as $m,n\to\infty$ to the entropy of $\mu$.
So it suffices to show the limit of the entropies is the entropy of the limit.
Let $W_k$ by a $k$-by-$k$ window in $\Z^2$ (or on $\Z^2_{m,n}$ when
$m$ and $n$ are larger than $k$).
Let $\mu_{m,n}(W_k)$ denote the restriction (projection)
of $\mu_{m,n}$ to $W_k$.
Then $\mu_{m,n}(W_k)$ converges to $\mu^*(W_k)$ by weak convergence. 
Since there are only a finite number of configurations on $W_k$,
$$\lim_{m,n\to\infty}H(\mu_{m,n}(W_k)) = H(\mu^*(W_k)).$$
Since $H(\mu_{m,n}(W_k))\geq H(\mu_{m,n})$
we have
$$H(\mu^*(W_k))\geq \lim_{m,n\to\infty}H(\mu_{m,n})=H_{top}.$$
Taking limits over $k$ gives $H(\mu^*)=H_{top}$ and therefore
by the theorem $\mu^*=\mu$.

We would like to thank Jeff Steif and Jim Propp for the idea
behind this argument. 

\subsubsection{The lozenge measure}
Burton and Pemantle's proof of Theorem \ref{unique} above uses a connection
between spanning trees and domino tilings. Lozenge tilings
are in fact also in bijection with a set of spanning trees, or rather,
directed spanning trees (``arborescences") on a directed triangular grid,
see \cite{PW}. However it is not evident how to 
adapt the proof of \cite{BP} to this case.
We will show nevertheless that the uniform measures on certain toroidal graphs
converge to a measure $\nu$, and their entropy converges
to the topological entropy. By the argument following Theorem \ref{unique},
the entropy of $\nu$ equals the topological entropy.
The only fact missing is the uniqueness of $\nu$.

\begin{figure}[htbp]
\PSbox{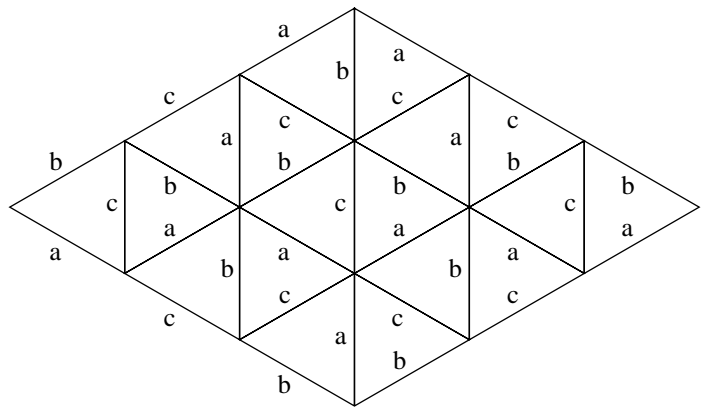}{.5in}{2.5in}
\caption{\label{33torus}The $18$ triangles make a torus.}
\end{figure}

Let $H_{m,n}$ be the quotient of $H$ by the lattice of translations
generated by $m\hat{x}$ and $n\hat{y}$. Then $H_{m,n}$ is a 
graph with $2mn$ vertices.
Let $\nu_{m,n}$ be the uniform measure on matchings of $H_{m,n}$.
\begin{lemma} The entropy of $\nu_{n,n}$ converges to the topological entropy
of lozenge tilings of the plane.
\label{loztori}
\end{lemma}
\begin{proof}
The proof is a simple variant of an unpublished trick due to G. Kuperberg.
Let $F_n$ be the set of free lozenge tilings of the region $R_n$, where $R_n$
consists of an equilateral triangle of side $n\sqrt{3}$
(recall that a lozenge has side length $\sqrt{3}$). 
Then $F_n$ is partitioned into equivalence classes,
two tilings being equivalent if they have the same set of tiles
which protrude beyond 
the boundary. The number of equivalence classes is less than $C^n$
for some constant $C$, so some equivalence class $B_n$ must satisfy
$$C^n|B_n|\geq |F_n|.$$

Reflect the region $R_n$ along one of its edges $c$.
Given any two tilings $T_1,T_2\in B_n,$ there is a tiling of the union of $R_n$ 
and its reflection, where the $R_n$ is filled with $T_1$ and its reflection
is filled with the reflection of $T_2$. 
This is because the tiles which protrude across the boundary edge $c$
are {\it fixed} under this reflection.

Take the union of $18$ reflections of $R_n$ so as to make
the parallelogram of Figure \ref{33torus}.
Given any $18$ elements of $B_n$, we can fill them in the $18$ triangles 
to make a tiling of this whole parallelogram.
Furthermore, the configurations of
tiles protruding beyond opposite boundaries of this parallelogram 
are translates of
each other. As a consequence we can glue opposite sides of the
parallelogram
together to make a tiling of the $3n$ by $3n$ torus (tilings of which 
correspond to perfect matchings of $H_{3n,3n}$).
We therefore have (where $N_{3n}$ is the set of tilings
of $H_{3n,3n}$)
$$(C^{-n}|F_n|)^{18}\leq  |B_n|^{18}\leq |N_{3n}|\leq |F_n|^{18}.$$
This implies that 
$$\lim_{n\to\infty}\frac1{n^2}\log|F_n|=\lim_{n\to\infty}\frac1{2n^2}
\log|N_{n}|,$$
which completes the proof.
\end{proof}

We will see later that $H(\nu_{m,n})$ and $H(\nu_{n,n})$ 
have the same limit as $m,n\to\infty$, 
so the restriction to the tori $H_{3n,3n}$ is unimportant. 
In fact we will show in section \ref{4.4} that there is a {\it unique}
weak limit to the $\nu_{m,n}$. Let $\nu$ be this limit.
We make the following
\medskip

\noindent{\bf Conjecture.}
{\it The measure $\nu$  is the {\rm unique}
measure of maximal entropy for lozenge tilings
of the plane. }

\subsection{Height functions.}
For any perfect matching of a simply connected subgraph $H'$ of $H$
one can associate an
integer-valued {\bf height function} to the faces of $H'$.
(By {\bf face} of $H'$ we mean a basic hexagon which contains an edge of $H'$.)

\begin{figure}[htbp]
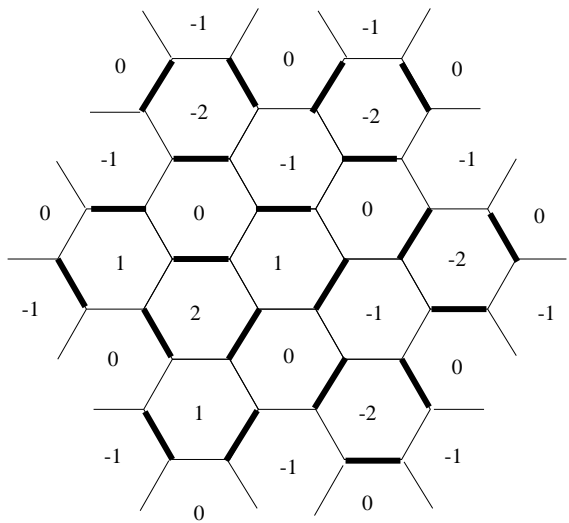

\PSbox{hts.eps}{1in}{3in}
\caption{\label{hts}Height function for lozenges.}
\end{figure}

Let $T$ be a perfect matching of $H'$.
Pick a face $x_0$ arbitrarily and assign it height $h(x_0)=0$.
For any other face $x$, 
take a path $x_0,x_1,\ldots,x_n=x$ of faces, where for each $j$,
$x_{j+1}$ is adjacent to $x_j$ along an unmatched edge.
Let $w=e^{2\pi i/3}$.
Then $h(x_{j+1})=h(x_j)\pm 1$, 
where $h$ increases by $1$ if $x_{j+1}-x_j$ is in directions 
$i,iw,$ or $iw^2$ and
$h$ decreases by $1$ if $x_{j+1}-x_j$ is in directions $-i,-iw$ or $-iw^2$.
This defines a height on each face, which can easily be shown to be
independent of the choice of path to that face. Figure \ref{hts}
shows a height function for the matching in Figure \ref{hexgraph}.

For dominos there is a similar definition; see \cite{Thu}. 

For a simply connected subgraph $H'$ which has a perfect matching, 
the height function along the faces bounding the region (that is, faces
outside $H'$ but containing an edge of $H'$) is well defined up to
an additive constant and is independent of the matching \cite{Thu}
(rather, this is true if the union of $H$ and these faces
is still simply connected).

\section{The measure of a cylinder set}
\label{cyls}
The results in this section are stated for lozenges, but they
apply verbatim for dominos, on condition that one puts weight $i=\sqrt{-1}$ on
vertical edges in the adjacency matrix $B$.

Let $H'$ be a finite, balanced,
simply connected subgraph of $H$,
with vertices $\{v_1,\ldots,v_n\}\in V_0$ and $\{v_1',\ldots,v_n'\}\in V_1$. 
Let $E=\{e_1,\ldots,e_k\}$ be a subset of disjoint edges of $H'$,
with $e_j=v_{p_j}v'_{q_j}$ for $j=1,\ldots,k$.

In the formula (\ref{detB}), each nonzero
term in the sum corresponds to a perfect
matching; those terms whose matchings contain all the edges in $E$ are
precisely the terms which contain
the subproduct 
$$a_E\stackrel{def}{=}a_{p_1q_1}a_{p_2q_2}\cdots a_{p_kq_k}.$$

The {\it sum} of the terms containing this product 
is (up to sign) the product of $a_E$ with 
the determinant of the cofactor $B_E$, that is,
the determinant of the matrix obtained from $B$ by
removing all rows $p_1,\ldots,p_k$ and all columns $q_1,\ldots,q_k$.
This can be proved inductively as follows. Expand $\det(B)$
along the row $p_1$:
$$\det(B)=(-1)^{p_1+1}(a_{p_11}\det B_{p_1,1}-a_{p_12}\det B_{p_1,2}+
\ldots\pm a_{p_1n}\det B_{p_1,n}).$$
The only term that contains $a_{p_1q_1}$ is the term for column $q_1$,
which is the term
$(-1)^{p_1+q_1}a_{p_1q_1}\det B_{p_1q_1}$. Now expand $B_{p_1q_1}$
along row $p_2$, and so on. Thus the sum of the terms containing $a_E$
is $(-1)^{\sum p_j+q_j}a_E\det(B_E)$.

Since $|a_E|=1$ this proves
\begin{prop} The number of perfect matchings containing all edges in $E$
is $$|\det(B_E)|.$$
\label{first}
\end{prop}

Later on the sign of $\det(B_E)$ will be important to us.
We note therefore that the number of perfect matchings containing
all the edges in $E$ is
$$(-1)^{\sum(p_j+q_j)}a_E\det(B_E)\mbox{sign}(\det B).$$

In a more general setting (e.g. in \cite{CKP})
one might wish to consider weighting the edges non-trivially, in which case
the formula $|a_E\det(B_E)|$
counts each matching according to the product of the weights
on its edges. 

Note that the quantity $\det(B_E)$
does not depend on the actual set of edges in $E$,
but only on the set of vertices involved. 
This shows, as expected, that each matching of the vertices
of $E$ can be extended to the same number of perfect matchings.

The matrix $B_E$ is an $(n-k)\times(n-k)$ cofactor of $B$; its 
determinant is equal to $\det(B)$ times
the determinant of the $k\times k$ cofactor of 
the inverse of $B$, $(B^{-1})_{E^*}$, 
where $E^*$ is the set of rows and columns not involved in $E$. 
So we have
\begin{thm} The number of perfect matchings containing all edges in $E$ is
$$|\det((B^{-1})_{E^*})\det(B)|.$$
That is, 
$$\mu_{H'}(U_E)=|\det((B^{-1})_{E^*})|.$$
More precisely, 
$$\mu_{H'}(U_E)=(-1)^{\sum p_j+q_j}a_E\det((B^{-1})_{E^*}).$$
\label{E}
\end{thm}

The advantage of this over Proposition \ref{first} is that the computation
is a determinant of size $k$ rather than $n-k$.
In particular if $E$ consists of a single edge $v_iv'_j$, 
the probability that that edge is in a random matching is
the absolute value of the $ij$-th entry of $B^{-1}$.

The rest of the paper consists of the computation of
$B^{-1}$ and applications of Theorem \ref{E}.
\section{Lozenge tilings}
\label{lozs}
\subsection{The torus}
Recall that  $H_{m,n}$ is the toroidal graph $H/(m\Z\hat{x}+n\Z\hat{y})$, 
and the vertex $r\hat{x}+s\hat{y}+t$ is designated by the triple 
$(r,s,t)\in\Z/m\Z\times\Z/n\Z\times\{0,1\}$.
Let $A_1$ be the adjacency matrix of $H_{m,n}$. 
Let $A_2$ be the matrix obtained from $A_1$ by changing the sign of the entries
corresponding to edges from $(m-1,k,1)$ to $(0,k,0)$ for each $k\in\Z/n\Z$.
Let $A_3$ be the matrix obtained from $A_1$ by changing the sign of the entries
corresponding to edges from $(k,n-1,1)$ to $(k,0,0)$ for each $k\in\Z/m\Z$.
Let $A_4$ be the matrix obtained from $A_1$ by changing the sign of 
both these sets of entries.

Let $B_1,B_2,B_3,B_4$ be the matrices obtained from $A_1,A_2,A_3,A_4,$ respectively,
as in (\ref{B}).

Following Kasteleyn, on condition that $m,n$ are both even, we have
\begin{lemma}[\cite{Kast1}]
\label{BBBB}
The number of perfect matchings $Z_{m,n}$ of $H_{m,n}$ is given by
\begin{equation}\label{BBBBeq}Z_{m,n}=\frac12(-\det B_1+\det B_2+\det B_3+\det B_4).
\end{equation}
\end{lemma}
A term for each perfect matching appears in each of these
four determinants.
The signs are arranged so that each perfect matching is counted
three times with sign $+$ and once with sign $-$ (see \cite{Kast1}).
For other parities of $m$ and $n$ a similar equation holds but with different
signs. Throughout the rest of this section we will assume
$m$ and $n$ are both even and neither is divisible by $3$.
This is simply a matter of convenience; the computations go through in the other
cases but require slightly different arguments.

Theorem \ref{E} also applies to tilings of the torus. The quantity 
$$\det B_1\det(B_1^{-1})_{E^*}$$
contains a term for each matching containing all edges of $E$.
The sign of these terms is the same as the sign of the corresponding
terms in $\det(B_1)$ (this easily follows from the proof). 
This is also true for $B_2,B_3,B_4$, and so
\begin{lemma} The number of perfect matchings of $T_{m,n}$ 
containing a given set of disjoint edges $E$ is given by 
\begin{eqnarray}
\label{probwith4}
\frac12\left|-\det B_1\det(B_1^{-1})_{E^*}+\det B_2\det(B_2^{-1})_{E^*}+
\right.\hspace{1in}\\
\hspace{1in}
+\left.\det B_3\det(B_3^{-1})_{E^*}+\det B_4\det(B_4^{-1})_{E^*}\right|.
\nonumber
\end{eqnarray}
\end{lemma}
\medskip

\subsection{Computation of the eigensystems.}
The eigenvectors of the $A_j$ are functions
$$f(r,s,t)=c_t z^{r}w^{s},~\mbox{and}~
f'(r,s,t)=c'_t z^rw^s$$
where 
\begin{eqnarray*}
c_0&=&\sqrt{1+z^{-1}+w^{-1}}\\
c_1&=&\sqrt{1+z+w}\\
c'_0&=&c_0\\
c'_1&=&-c_1
\end{eqnarray*}
(for either choice of sign of the square root), and
$$
\begin{array}{lll}
z^m=1,&w^n=1&\mbox{for}~A_1\\
z^m=-1,&w^n=1&\mbox{for}~A_2\\
z^m=1,&w^n=-1&\mbox{for}~A_3\\
z^m=-1,&w^n=-1&\mbox{for}~A_4.
\end{array}
$$
The corresponding eigenvalues are 
\begin{eqnarray*}
\lambda&=&c_0c_1\\
\lambda'&=&-c_0c_1.
\end{eqnarray*}

Consider for example the matrix $A_1$.
Let 
\begin{eqnarray*}
z_1&=&e^{2\pi i/m}\\
w_1&=&e^{2\pi i/n}\\
c_{k,\ell,0}&=&\sqrt{1+z_1^{-k}+w_1^{-\ell}}\\
c_{k,\ell,1}&=&\sqrt{1+z_1^k+w_1^\ell}\\
c'_{k,\ell,0}&=&c_{k,\ell,0}\\
c_{k,\ell,1}'&=&-c_{k,\ell,1}
\end{eqnarray*}
(any choice of sign for the square roots will do).
Then $f_{k,\ell}(r,s,t)=c_{k,\ell,t}z_1^{kr}w_1^{\ell s}$ and
$f'_{k,\ell}(r,s,t)=c'_{k,\ell,t}z_1^{kr}w_1^{\ell s}$ are eigenfunctions.
The functions $f_{k,\ell}+f'_{k,\ell}$ are zero on $V_1$,
and run through a basis for all functions on $V_0$ (the exponentials)
as $(k,\ell)$ runs through $\{0,\ldots,m-1\}\times\{0,\ldots,n-1\}$.
Similarly the functions $f_{k,\ell}-f'_{k,\ell}$
are zero on $V_0$ and run through a basis for functions on $V_1$.
Thus the collection of $f_{k,\ell}$ and $f'_{k,\ell}$ is a
basis of eigenvectors for the matrices $A_1$.
In particular we have
\begin{eqnarray}
\nonumber\det(A_1)&=&\prod_{k=0}^{m-1}\prod_{\ell=0}^{n-1}\lambda_{k,\ell}
\lambda'_{k,\ell}\\
&=&\prod_{k=0}^{m-1}\prod_{\ell=0}^{n-1}\nonumber
-(1+e^{2\pi ik/m}+e^{2\pi i\ell/n})(1+e^{-2\pi ik/m}+e^{-2\pi i\ell/n}).\\
&=&\prod_{k=0}^{m-1}\prod_{\ell=0}^{n-1}\label{detA1}
(1+e^{2\pi ik/m}+e^{2\pi i\ell/n})^2
\end{eqnarray}
(we removed the initial minus sign in this last equality
since there are an even number of terms,
and also in the second product we replaced $k,\ell$ with $-k,-\ell$).
Similarly one can compute
\begin{eqnarray}
\nonumber\det(A_2)&=&\prod_{k,\ell}
(1+e^{\pi i(2k+1)/m}+e^{2\pi i\ell/n})
(1+e^{-\pi i(2k+1)/m}+e^{-2\pi i\ell/n})\\ \label{A2}
&=&\prod_{k,\ell}
(1+e^{\pi i(2k+1)/m}+e^{2\pi i\ell/n})^2
\end{eqnarray}
and furthermore 
\begin{eqnarray}
\label{A3}
\det(A_3)&=&\prod_{k,\ell}
(1+e^{2\pi ik/m}+e^{\pi i(2\ell+1)/n})^2\\
\label{A4}
\det(A_4)&=&\prod_{k,\ell}
(1+e^{\pi i(2k+1)/m}+e^{\pi i(2\ell+1)/n})^2.
\end{eqnarray}
These four determinants
are non-zero by our assumption that
$m$ and $n$ are nonzero modulo $3$ (three complex numbers of modulus $1$
sum to zero only if they are a rotation of the set 
$\{1,e^{2\pi i/3},e^{4\pi i/3}\}$).

\subsection{Computation of the inverses of $B_1,\ldots,B_4$.}
It suffices to invert the $A_j$.
Define the vector
$$\delta_{0,0,0}(x,y,t) 
=\left\{\begin{array}{ll}1&\mbox{if}~(x,y,t)=(0,0,0)\\
0&\mbox{otherwise.}\end{array}\right.$$
Consider first the matrix $A_1$.
We can write
$$\delta_{0,0,0}(x,y,t)= \frac{1}{2mn}
\sum_{k,\ell}\frac{f_{k,\ell}(x,y,t)+f_{k,\ell}'(x,y,t)}{c_{k,\ell,0}}.$$
Thus the vector 
\begin{equation}\label{P00}
P_{0,0,0}(x,y,t)=\frac{1}{2mn}\sum_{k,\ell}\frac1{c_{k,\ell,0}}
\left(\frac{f_{k,\ell}(x,y,t)}{\lambda_{k,\ell}}+\frac{f_{k,\ell}'(x,y,t)}
{\lambda'_{k,\ell}}\right)
\end{equation}
has the property that 
$$A_1P_{0,0,0}(x,y,t)=\delta_{0,0,0}(x,y,t).$$
The preimages of delta functions at other vertices can be similarly defined: 
define
$P_{r,s,0}(x,y,t)=P_{0,0,0}(x-r,y-s,t)$;
then $$A_1P_{r,s,0}(x,y,t)=\delta_{r,s,0}(x,y,t)=
\left\{\begin{array}{ll}1&\mbox{if}~(x,y,t)=(r,s,0)\\
0&\mbox{otherwise.}\end{array}\right.$$
Thus the function $P_{r,s,0}$ gives the entries of the
$(r,s,0)-$column of the matrix $A_1^{-1}$.

For notational simplicity we will write $P^{(1)}(x,y,t)=P_{0,0,0}(x,y,t)$
(the superscript $(1)$ refers to the fact that it comes from
matrix $A_1$; note that $P^{(1)}$ also depends on $m$ and $n$).

If we simplify (\ref{P00}) above:
$$P^{(1)}(x,y,t)=\frac{1}{2mn}\sum_{k,\ell}\frac1{c_{k,\ell,0}}
\left(\frac{f_{k,\ell}(x,y,t)-f_{k,\ell}'(x,y,t)}{\lambda_{k,\ell}}\right),$$
we see that $P^{(1)}(x,y,0)=0$ for all $x,y$.
On the other hand plugging in for $f,c,\lambda,z_1,w_1$ we have
\begin{equation}\label{Psum}
P^{(1)}(x,y,1)=\frac{1}{mn}\sum_{k,\ell}\frac{e^{2\pi ikx/m}e^{2\pi i\ell y/n}}{
1+e^{-2\pi ik/m}+ e^{-2\pi i\ell/n}}.
\end{equation}

Similar expressions hold for the inverses of $A_2,A_3,A_4$, except that the
sum is over a different set of angles.
For example for $A_2$ we have 
$$P^{(2)}(x,y,1)= \frac{1}{mn}\sum_{k,\ell}\frac{e^{\pi i(2k+1)x/m}
e^{2\pi i\ell y/n}}{ 1+e^{-\pi i(2k+1)/m}+ e^{-2\pi i\ell/n}}.$$

\subsection{The limiting measure.}
Here we compute the 
limit of the ratio of (\ref{probwith4}) and (\ref{BBBBeq}), that is,
the limit of the cylinder set measures $\nu_{m,n}(U_T)$.
We still assume $m,n\not\equiv 0\bmod 3$ so that the determinants in
(\ref{detA1})-(\ref{A4}) are nonzero. 

We will show below that the quantities 
$\det((B_j^{-1})_{E^*})$ for $j=1,2,3,4$ 
all converge to the same value as $m,n\to\infty$.
Then (if we remove the absolute value sign in (\ref{probwith4}))
the ratio of (\ref{probwith4}) and
(\ref{BBBBeq}) is a weighted average of these four quantities, with weights
$\pm\det B_j/Z_{m,n}$.
These weights are all in the interval $[-1,1]$ since
$Z_{m,n}\geq |\det B_j|$ for each $j$: recall that $Z_{m,n}$ and $|\det{B_j}|$ count
the same objects except that some signs in $|\det{B_j}|$ are negative.
Since the weights sum to $1$,
the weighted average is also converges to the same value as each
$\det((B_j^{-1})_{E^*})$. 

It suffices to show that, for fixed $x,y$, for each $j=1,2,3,4$
the quantities $P^{(j)}(x,y,1)$ tend to the same
value (recall that the matrix $B_j^{-1}$ 
has entries which are the values of $P^{(j)}$).

For a fixed integers $x,y$, the function 
$$
p(\theta,\phi)=\frac{e^{ix\theta}e^{iy\phi}} {1+e^{-i\theta}+e^{-i\phi}}$$
is continuous on $([0,2\pi]\times [0,2\pi])\setminus
\{(2\pi/3,4\pi/3)\cup(4\pi/3,2\pi/3)\}$,
and has a pole at the two points $(2\pi/3,4\pi/3)$ and $(4\pi/3,2\pi/3)$.
The expression (\ref{Psum}) and the corresponding expressions for
the other $P^{(j)}$ are each Riemann sums for the integral
of $p$ over $[0,2\pi]\times[0,2\pi]$. 
Fix $\delta>0$ and let $N_\delta$ be a 
$\delta\times\delta$-neighborhood of the pole $(\theta_0,\phi_0)=(2\pi/3,4\pi/3)$.
Outside $N_\delta$ (and the corresponding neighborhood of the other pole)
the Riemann sums converge to
the appropriate integral. We must show that on $N_\delta$
the sums are small (in the sense that they tends to zero with $\delta$ 
as $m,n\to\infty$).

On $N_\delta$, we use the Taylor expansion of the denominator:
$$
\frac1{mn}\left|\sum_{(\theta,\phi)\in N_\delta}p(\theta,\phi)\right|\leq
\hspace{2in}$$
\begin{eqnarray*}
\leq &
\frac1{mn}\sum_{N_\delta}\frac{1}{|-ie^{-i\theta_0}(\theta-\theta_0)-
ie^{-i\phi_0}
(\phi-\phi_0)+O((\theta-\theta_0)^2,(\phi-\phi_0)^2)|}\\
=&\frac1{mn}\sum_{N_\delta}
\frac1{|e^{-2\pi i/3}(\theta-\theta_0)+e^{-4\pi i/3}(\phi-\phi_0)+
O((\theta-\theta_0)^2,(\phi-\phi_0)^2)|}\\
=&\frac1{mn}\left[\sum_{N_\delta}
\frac1{|e^{-2\pi i/3}(\theta-\theta_0)+e^{-4\pi i/3}(\phi-\phi_0)|}\right]
+O(\theta-\theta_0,\phi-\phi_0).
\end{eqnarray*}
(The reason we can take the big-$O$ term out of the denominator
is that 
$$|e^{-2\pi i/3}(\theta-\theta_0)+e^{-4\pi i/3}(\phi-\phi_0)|\geq
const(|\theta-\theta_0|+|\phi-\phi_0|)$$
for $|\theta-\theta_0|$ and $|\phi-\phi_0|$ small; we can then
move it out of the summation altogether because of the 
factor $\frac1{mn}$.)

Since $m,n\not\equiv 0\bmod 3$, we have, in the sum for each $P^{(j)}$,
that $\theta-\theta_0$ is of the form
$\theta-\theta_0=\pi(\frac km-\frac23)
=\pi\frac{k'}{3m}$ and similarly $\phi-\phi_0=\frac{\pi \ell'}{3n}$,
where $k'$ and $\ell'$ are not multiples of $3$.
Thus the sums are each bounded by one of the form
$$\frac1{mn}\sum
\frac{1}{|\frac k{3m}+\frac \ell{3n}e^{4\pi i/3}|} + O(\delta),$$
where the sum is over 
$$\left\{(k,\ell)\in\Z^2\mid~(\frac k{3m},\frac \ell{3n})\in 
N_\delta\setminus\{(0,0)\}\right\}.$$
This sum is easily shown to be
bounded independently of $m,n$ and tending to zero as $\delta\to 0$:
this follows from the integrability near the origin of the function 
$\frac1{|x+e^{4\pi i/3}y|}$ with respect to $dxdy$.

\begin{thm}\label{lozcouple}
For a finite matching $E$ 
covering black vertices $\{b_1,\ldots,b_k\}$ and white vertices
$\{w_1,\ldots,w_k\}$, we have 
$\nu(U_E)=|\det P(w_i-b_j)|$, where 
$$P(x,y,1)=\frac{1}{4\pi^2}\int_0^{2\pi}\int_0^{2\pi}
\frac{e^{ix\theta}e^{iy\phi}d\theta d\phi} {1+e^{-i\theta}+e^{-i\phi}}.$$
\end{thm}

\begin{proof}
This follows from Theorem \ref{E} and the convergence argument above.
\end{proof}

\label{4.4}
\subsection{The coupling function for lozenge tilings of the plane.}
We call $P(x,y,t)$ the {\bf coupling function} for lozenge tilings of the plane.
Let $z=e^{-i\theta}$ and $w=e^{-i\phi}$.
Then we obtain the contour integral
\begin{equation}\label{Peq}
P(x,y,1)=\frac{1}{4\pi^2}\int_{S^1\times S^1}
\frac{z^{-x}w^{-y}}{(1+z+w)}\frac{dz}{(-iz)}\frac{dw}{(-iw)}.
\end{equation}

The function $P$ has all the symmetries of the graph $H$, that is,
\begin{equation}\begin{array}{ccc}\label{symmetry}
P(x,y,1)&=&P(-x-y-1,x,1)=P(y,-x-y-1,1)\\
=P(y,x,1)&=&P(-x-y-1,y,1)=P(x,-x-y-1,1).
\end{array}
\end{equation}
These are obtained by noticing that (\ref{Peq}) is invariant
under interchanging $z$ and $w$,
or under the substitution $(z,w)\to(w^{-1},zw^{-1})$, or under combinations of
these.

We can explicitly evaluate the integral (\ref{Peq}) as follows.
This is slightly easier to calculate in the case $x\leq -1$. 
The other values can then be obtained by (\ref{symmetry}).
If we fix $w$ and integrate over $z$, there is a pole inside 
the unit circle when $|1+w|<1$, that is, when $-\phi\in(2\pi/3,4\pi/3)$.
The residue of $z^{-x-1}/(1+z+w)$ at this pole is $(-1-w)^{-x-1}$.
For $-\phi\not\in(2\pi/3,4\pi/3)$ there is no pole 
inside the circle and therefore the integral is zero.

\begin{figure}[htbp]
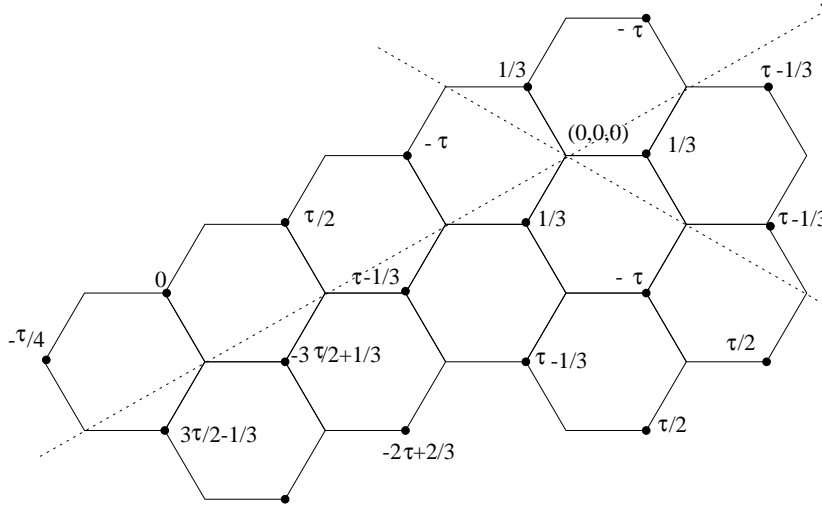

\PSbox{hexvals.eps}{0in}{3in}
\caption{\label{hexvals}Values of $P(x,y,t)$. Here $\tau=\sqrt{3}/(2\pi)$.}
\end{figure}

So when $x\leq -1$, 
\begin{equation}
\label{pint}
P(x,y,1)=\frac{-i}{2\pi}(-1)^{-x-1}
\int_{e^{2\pi i/3}}^{e^{4\pi i/3}}w^{-y-1}(1+w)^{-x-1} dw.
\end{equation}
For example when $x=-1$ we have
\begin{equation}\label{p1n0}
P(-1,y,1)=\frac{i}{2\pi}\left(\frac{e^{2\pi iy/3}-e^{4\pi iy/3}}{y}\right)=
\frac{c_y\sqrt{3}}{2\pi y},
\end{equation}
where $c_y=0,1,-1$ if $y\equiv 0,-1,1 \bmod 3$.
Some values of $P(x,y,1)$ are plotted in Figure \ref{hexvals}.
To compute these values, start from (\ref{p1n0}), which, along with the symmetries
(\ref{symmetry}) gives the value on $3$ lines, and use the fact that
the sum of the three values adjacent to any vertex 
(except the origin) is zero (since $AP=0$ at these vertices).
For example the equation 
$$P(-k,0,1)+P(-k+1,-1,1)+P(-k+1,0,1)=0$$ 
determines $P(-k,0,1)$ for $k>0$ 
inductively for increasing $k$ 
(note as before that $P(-k,-1,1)=P(-1,-k,1)=\frac{c_k\sqrt{3}}{2\pi k}$).

\subsection{Some local probabilities for lozenges}

From Theorem \ref{E} and equation (\ref{pint}) with $x=-1,y=0$,
the probability of the edge $(0,0,0)(0,0,1)$
being in a random matching is given by $P(0,0,1)=\frac13$, as expected
by symmetry.
The probability of both the edges $(0,0,0)(0,0,1)$ and
$(0,1,0)(0,1,1)$ occurring is given by 
$$\left|\begin{array}{cc}P_{0,0,0}(0,0,1)&P_{0,0,0}(0,1,1)\\
P_{0,1,0}(0,0,1)&P_{0,1,0}(0,1,1)\end{array}\right|
=\left|\begin{array}{cc}\frac13&\frac{\sqrt{3}}{2\pi}-\frac13\\
\frac13&\frac13\end{array}\right|=
\frac29-\frac{\sqrt{3}}{6\pi}\approx 13\%.$$

The probability of the three edges 
$$(0,0,0)(0,0,1),~~~ (0,1,0)(-1,1,1)\mbox{ and }(-1,1,0)(-1,0,1)$$
(which bound a hexagon) all occurring, is:
$$\left|\begin{array}{ccc}-\frac{\sqrt{3}}{2\pi}&\frac13&\frac13\\
\frac13&-\frac{\sqrt{3}}{2\pi}&\frac13\\
\frac13&\frac13&-\frac{\sqrt{3}}{2\pi}\end{array}\right|=
\frac{2}{27}+\frac{\sqrt{3}}{6\pi}
-\frac{3\sqrt{3}}{8\pi^3}\approx 14.5\%.$$
(Such a hexagon corresponds to a local minimum for the height function.)

Note from (\ref{p1n0}) that $P(-3k,3k,1)=P(-1,-3k,1)=0$.
As a consequence, the set of edges $E=\{e_n\}_{n\in\Z}$ where $e_n$ is the edge
$(-3n,3n,0)(-3n,3n,1)$ is an ``independent" set: the probability of any
fixed subset of size $k$
occurring in a random matching is exactly $(1/3)^k$.
Is there any simple explanation for this fact?
\medskip

We can also compute long-range correlations.
The probability of the horizontal edge $(0,0,0)(0,0,1)$ and the horizontal
edge $(-n,n,0)(-n,n,1)$ both occurring is given by the determinant
$$\left|\begin{array}{cc}\frac13&P(-n,n,1)\\
P(n,-n,1)&\frac13 \end{array}\right|=\frac 19-P(-n,n,1)^2.$$
By (\ref{symmetry}) and (\ref{p1n0}), 
$P(-n,n,1)=P(-1,-n,1)=c_n\sqrt{3}/(2\pi n).$
So the probability of these two edges both occurring
tends quadratically (but not faster) to $\frac19$.
More generally we have
\begin{prop} \label{1/n}
$$P(x,y,1)=O\left(\frac1{|x|+|y|}\right).$$
\end{prop}

\begin{proof} We show for $x\leq -1$ that
$P(x,y,t)=O(\frac1{|x|})$; by the dihedral $6$-fold
symmetry of $P$ the result will then follow. From (\ref{pint}) we have
\begin{eqnarray*}
|P(x,y,1)|&\leq&\frac1{2\pi}\int_{e^{2\pi i/3}}^{e^{4\pi i/3}}
|1+w|^{-x-1}|dw|\\
&=&\frac2{2\pi}\int_{2\pi/3}^\pi(2\cos(\theta/2))^{-x-1}d\theta\\
&=&\frac1\pi\int_0^1\frac{t^{-x-1}}{\sqrt{1-(t/2)^2}}dt\\
&\leq& \frac2{\pi}\int_0^1 t^{-x-1}dt = O(\frac1{|x|})
\end{eqnarray*}
where we used the substitution $t=2\cos(\theta/2)$.
\end{proof}

Let $E_1$ and $E_2$ be two finite collections of edges. The probability that
edges $E_1\cup E_2$ are present in a matching is given by a determinant
\begin{equation}\label{detexp}
\left|\begin{array}{cc}C_1&D_1\\D_2&C_2\end{array}\right|,
\end{equation}
where $\det(C_1)=\mu(U_{E_1})$ and $\det(C_2)=\mu(U_{E_2})$.
The entries of $D_1$ and $D_2$ are the values of $P$ connecting
points of $E_1$ to $E_2$. If the distance between $E_1$ and $E_2$
is $\geq n$ then these entries are all $O(\frac1n)$. 
Let $\sum \epsilon(\sigma)a_{1\sigma(1)}\ldots a_{k\sigma(k)}$
be the expansion of the determinant (\ref{detexp}).
Suppose $a_{i\sigma(i)}$ is an entry from the submatrix $D_1$, that
is, $i\leq|E_1|$ and $\sigma(i)>|E_1|$. Then there must exist
an index $j$ with $j>|E_1|$ and $\sigma(j)\leq |E_1|$ (since 
$\sigma$ is a bijection), that is,
$a_{j\sigma(j)}$ is an entry of $D_2$. We may similarly conclude
that each term in the sum for the determinant
has the same number of entries from $D_1$ as from $D_2$.
Therefore the determinant (\ref{detexp})
is equal to $\det(C_1)\det(C_2)+O(\frac1{n^2})$.

This proves 
\begin{thm}\label{mixing}If the distance between $E_1$ and $E_2$
is at least $n$, then
$$\mu(U_{E_1}\cap U_{E_2})=\mu(U_{E_1})\mu(U_{E_2})+O(\frac1{n^2}).$$
\end{thm}

Note that the comments before Proposition (\ref{1/n}) 
give a corresponding quadratic lower bound if $E_1$ and $E_2$ are
single horizontal tiles separated by $-(3n+1)\hat{x}+(3n+1)\hat{y}$.

An argument similar to Theorem \ref{mixing}
shows that $\nu$ is {\bf mixing of all orders},
that is, given $k$ events $E_1,\ldots,E_k$ each of which is
separated from the others, the joint probability converges
to the product of the individual probabilities as the 
distance between them tends to infinity.

\subsection{The height distribution}

Let $h_n$ be the random variable which gives the height 
difference between the face at the origin 
(by which we mean, the face above and right of the point $(0,0,0)$) 
and the 
face at the point $(-n,n,0)$ in a random lozenge 
tiling of the plane. In this section we compute
the variance of $h_n$.

Since the height is only defined up to an additive constant,
we may as well assume that the height on the face at the origin is zero.

One motivation for computing this variance is that it is related to the
number of cycles in the union of two
random matchings as follows.
Given two random matchings $M_1$ and $M_2$, 
their union is a set of disjoint cycles.
The difference of their height functions $h^{(1)}-h^{(2)}$ has the property
that it is constant on the connected
components of the complement of these cycles,
and changes by $3$ when crossing any
of these cycles. The height difference across a cycle
increases or decreases by $3$ depending on which
of the two possible matchings of the cycle occurs in $M_1$ (with 
the complementary matching necessarily occurring in $M_2$).
Since the two matchings were chosen randomly, 
across each cycle the difference of heights increases or
decreases by $3$ {\it independently} of what happens at the other cycles.
(Among all pairs of matchings whose union has this particular cycle structure,
each matching is equally likely to have either ``half" of each cycle.) 

Thus for a face $f$ separated from the face at the origin by $k$ cycles, 
the variance in the height difference 
(that is, the variance $\sigma^2(\Delta h(f))=\sigma^2(h^{(1)}(f)-h^{(2)}(f))$)
is the same as the variance of a sum of $k$ independent fair 
coin flips of values $\pm 3$.  This variance is $9k$. 
Since $M_1$ and $M_2$ were chosen randomly, the number $y$ of cycles
separating $f$ from the origin is also a random variable.
We have $\sigma^2(\Delta h|y)=9y$,
and 
$$\sigma^2(\Delta h|y)=\E((\Delta h)^2|y)-\E(\Delta h|y)^2=\E((\Delta h)^2|y)$$
since $\E(\Delta h|y)=0$. Furthermore,
$$\sigma^2(\Delta h)=\E((\Delta h)^2)=\sum_y \P(y)\E((\Delta h)^2|y)
=\sum_y\P(y)9y = 9\E(y).$$
Thus in general the variance of 
the height differences at a fixed face (for two matchings)
is $9$ times the expected number of cycles separating that face from the origin.

Since the two matchings $M_1,M_2$ are chosen independently, the variance
of their 
height difference at $f$ is exactly twice the variance of the height
of a single matching at $f$
(note that the height at $f$ has mean value zero,
which follows from the fact that 
the probability of any edge being present is $1/3$).
We conclude
\begin{prop}\label{contours}
In the union of two random matchings, the expected number of 
cycles separating two faces $f_1,f_2$ is $2s/9$,
where $s$ is the variance in the height (between $f_1$ and $f_2$)
of a single random matching. 
\end{prop}

Let $E_n$ denote the set of horizontal edges $(-k,k,0)(-k,k,1)$
for $k=1,\ldots,n$.
\begin{lemma}\label{height} The value of $h_n$ 
is equal to $n-3r_n$, where $r_n$ is the number
of horizontal edges from $E_n$ in the matching. 
\end{lemma}

\begin{proof} This follows from the definition of the height function.
On the vertical path of faces from $(0,0,0)$ to $(-n,n,0)$, the height
increases by $1$ for edge of $E_n$ not present in the matching,
and decreases by $2$ for every edge present in the matching.
So the height difference is $(n-r_n)-2r_n$.
\end{proof}
\medskip

Let $M=M_n$ be the $n\times n$ 
matrix $M=(m_{k,j})$ where 
$$m_{k,j}=P(-|k-j|,|k-j|,1).$$
Then as we saw, the determinant of $M$ is the probability
that all the edges from $E_n$ are present in a random matching. 
By symmetry (\ref{symmetry}), when $\ell>0$, 
$$P(-\ell,\ell,1)=P(-1,-\ell,1)=\frac{c_\ell \sqrt{3}}{2\pi\ell}.$$
For example when $n=5$ we have (with $\tau=\frac{\sqrt{3}}{2\pi}$)
$$M_5=\left(\begin{array}{ccccc}\frac13&-\tau&\frac \tau2&0&\frac{-\tau}4\\
-\tau&\frac13&-\tau&\frac \tau2&0\\
\frac \tau2&-\tau&\frac13&-\tau&\frac \tau2\\
0&\frac \tau2&-\tau&\frac13&-\tau\\
\frac{-\tau}4&0&\frac \tau2&-\tau&\frac13
\end{array}\right).$$

For each subset $S\subset E_n$ of cardinality $k$, 
let $M_S$ be the submatrix of $M$ whose rows and columns are
those of $M$ indexed by $S$. Then $(-1)^{\sum p_j+q_j}\det(M_S)$ 
is the probability
that the edges corresponding to $S$ are present in a random tiling,
where the sum is over $S$.
Since $p_j=q_j$ in our ordering, the probability in question
is exactly $\det(M_S)$.

Let 
$$p(z)=z^n-\alpha_1z^{n-1}+\alpha_2z^{n-2}-\ldots+(-1)^n\alpha_n$$
be the characteristic polynomial
of $M$. Then $\alpha_k$ is the sum over all subsets $S\subset E_n$ of size
$k$ of $\det M_S$, that is, the sum over $S$ of
the probability that $S$ is in the random matching.

Let $q(z)=\sum_{k=0}^n \beta_kz^k$ where $\beta_k$ is the probability that there
are exactly $k$ edges from $E_n$ present.
A formula of Ch. Jordan \cite{Comtet} relates the $\alpha_k$ to the $\beta_k$:
$$\beta_k=\alpha_k-\left(\begin{array}{c}k+1\\k\end{array}\right)\alpha_{k+1}
+\left(\begin{array}{c}k+2\\k\end{array}\right)\alpha_{k+2}-\ldots.$$
From this we can easily derive the following:
\begin{thm} We have
$$q(z)=(1-z)^np(\frac1{1-z}).$$
\hfill{$\square$}
\end{thm}

Now if $p(z)=\prod_{i=1}^n(z-\lambda_i)$, then
$$q(z)=\prod_{i=1}^n(1-\lambda_i+\lambda_iz).$$

The expected value of $r_n$ is $\sum k\beta_k=q'(1)$.
We find
$$q'(z)=q(z)\sum_{i=1}^n\frac{\lambda_i}{1-\lambda_i+\lambda_iz},$$
so 
$$q'(1)=q(1)\sum\lambda_i=\sum\lambda_i=\mbox{trace}(M)=n/3.$$
This is as expected since the probability of each edge $(-k,k,0)(-k,k,1)$
being present is $1/3$. 

Similarly, the variance of $r_n$ is given by
$$\sigma^2(r_n)=
\sum k^2\beta_k - (\sum k\beta_k)^2=q'(1)+q''(1) - (\frac n3)^2.$$
We have
$$q''(z)=q(z)\sum_{i\neq j}\frac{\lambda_i\lambda_j}{(1-\lambda_i+\lambda_iz)
(1-\lambda_j+\lambda_jz)},$$
and so 
$$q''(1)=\sum_{i\neq j}\lambda_i\lambda_j=\mbox{trace}(M)^2-\mbox{trace}(M^2)=
(n/3)^2-\mbox{trace}(M^2).$$

To compute $\mbox{trace}(M^2)$, note that $M$ is symmetric and the
$i$th row for $1<i<n$ is:
$$(\frac{c_i\tau}{i-1},
\ldots,\frac \tau2,-\tau,\frac13,-\tau,\frac \tau2,\ldots,\frac{c_{n-i}\tau}{n-i}).$$
The inner product of this row with itself is
\begin{eqnarray*}
\frac19&+&\tau^2(1+\frac1{2^2}+\frac1{4^2}+\frac{1}{5^2}+\ldots+
\frac{c_{n-i}^2}{(n-i)^2})\\
&&+\tau^2(1+\frac1{2^2}+\frac1{4^2}+\frac{1}{5^2}+\ldots+
\frac{c_{i}^2}{(i-1)^2}).
\end{eqnarray*}
Similar sums (with $\frac{c_0\tau}{0}$ 
replaced with $1/3$) hold for the rows $i=1$ and $i=n$.
Summing these up for $i=1,\ldots,n$ yields
$$
\mbox{trace}(M^2)=
\frac n9+\tau^2\left(2(n-1)+2\frac{(n-2)}{2^2}+2\frac{(n-4)}{4^2}+\ldots+
2\frac{c_{n-1}^2}{(n-1)^2}\right)
$$
\begin{equation}
=\frac n9+\frac{3}{4\pi^2}\left(2n(1+\frac1{2^2}+\frac1{4^2}+\ldots+
\frac{|c_{n-1}|}{(n-1)^2})-
2(1+\frac12+\frac14+\ldots+\frac{|c_{n-1}|}{n-1})\right).
\label{sss}
\end{equation}
Note that 
$$
1+\frac1{2^2}+\frac1{4^2}+\ldots+\frac{|c_r|}{r^2}=
$$
\begin{eqnarray*}
&=&\left(1+\frac1{2^2}+\frac1{3^2}+\ldots+ \frac1{r^2}\right)-
\frac19\left(1+\frac1{2^2}+\frac1{3^2}+\ldots+
\frac{1}{\lfloor r/3\rfloor^2}\right)\\
&=&(1-\frac19)(\frac{\pi^2}{6}+O(\frac1r))
\end{eqnarray*}
and similarly
$$
1+\frac12+\frac1{4}+\ldots+\frac{|c_r|}{r}=
$$
\begin{eqnarray*}
&=&\left(1+\frac1{2}+\frac1{3}+\ldots+\frac1{r}\right)-
\frac13\left(1+\frac1{2}+\frac1{3}+\ldots+\frac1{\lfloor r/3\rfloor}\right)\\
&=&\log(r)+O(1)-\frac13(\log(r/3)+O(1))\\
&=&\frac23\log(r)+O(1).
\end{eqnarray*}
Plugging this into (\ref{sss}) gives
$$\mbox{trace}(M^2)= \frac n3-\frac{\log(n)}{\pi^2}+O(1),
$$
which yields 
$$\sigma^2(r_n)=\frac{\log(n)}{\pi^2}+O(1).$$
Finally we find
\begin{thm}\label{var}
Let $h_n$ be the height difference between the origin and the point
$(-n,n,0)$. The expected value of $h_n$ is $0$.
The variance in $h_n$ is given by $\sigma^2(h_n)=\frac{9\log(n)}{\pi^2}+O(1)$.
\end{thm}

The only previously known result in this vein is in \cite{CEP}, 
where Cohn, Elkies, and Propp found a $O(\sqrt{n})$ upper bound
for the height variance on a bounded region.

\section{Computation of the coupling function for dominos}
\label{dominopf}
\begin{proof}[Proof of Theorem \ref{domc}]
The proof follows from Theorem \ref{E},
the fact that $\mu_{m,n}\to\mu$ (see the comments
after Theorem \ref{unique}), and a convergence argument analogous
to that preceding Theorem \ref{lozcouple}.
\end{proof}

In fact a computation analogous to the case of lozenges yields:
\begin{thm}\label{Pdoms} The coupling function for dominos is
$$P(x,y)=
\frac{1}{4\pi^2}\int_0^{2\pi}\int_0^{2\pi}\frac{e^{i(x\theta+y\phi)}d\theta
d\phi}{2\cos(\theta)+2i\cos(\phi)}.$$
\end{thm}
This statement depends on our choice of weighting the vertical edges of 
$\Z^2$ with weight $i$.

\begin{figure}[htbp]
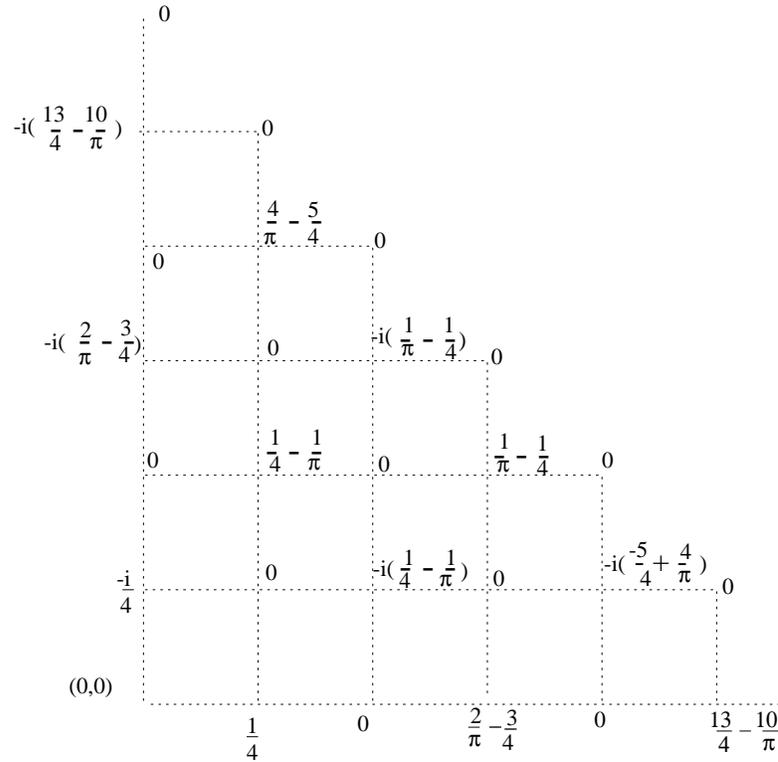

\PSbox{domvals.eps}{.5in}{4in}
\caption{\label{domvals}Values of $P(x,y)$ for dominos.}
\end{figure}

As for lozenges, 
$P(x,y)$ is zero for black vertices (i.e. when $x+y\equiv 0\bmod 2$).
Some values are shown in Figure \ref{domvals} (where the origin is
in the lower left corner).

One can explicitly evaluate, for $x\geq 1,$
$$P(2x+1,2x)= 
(-1)^x\left[\frac14-\frac1{\pi}\left(1-\frac13+\frac15-\ldots+(-1)^{x+1}\frac1{2x-1}
\right)\right],$$
and
$$P(2x,2x-1)=-iP(2x+1,2x).$$
From these values (and the diagonal symmetry)
the others can be obtained inductively using $AP=0$ except at the origin.

As an example, the probability of the two edges $(0,0)(1,0)$
and $(0,1)(1,1)$ both occurring is 
$$\left|\begin{array}{cc}\frac14&\frac{-i}4\\ \frac{-i}4&\frac14\end{array}
\right|=\frac 18.$$

The probability of the two edges $(0,0)(1,0)$ and $(0,1)(0,2)$
both occurring is given by
$$i\cdot\left|\begin{array}{cc}\frac14&\frac{-1}{\pi}+\frac 14\\
\frac{-i}4&\frac{-i}4\end{array} 
\right|=\frac{1}{4\pi}.
$$ 
This probability was also computed in \cite{BP}, although their stated
value is incorrect (they have acknowledged 
the mistake in their arithmetic).

Proposition \ref{1/n} and Theorem \ref{mixing} also apply to dominos.
However the computation of the variance in the height function
cannot be done within this framework since there is no simple analogue
of Lemma \ref{height}.

\bigskip
\noindent{\it Acknowledgement.} I 
would like to thank Oded Schramm, Jeff Steif, Jim Propp
and the referee for many helpful suggestions regarding
this paper.

\end{document}